\documentstyle[12pt]{article}
\begin{document}
\newtheorem{theorem}{Theorem}
\newtheorem{lemma}{Lemma}
\newtheorem{definition}{Definition}
\newtheorem{proposition}{Proposition}
\title{Plemelj Projection Operators over Domain Manifolds}
\author{John Ryan\\Department of Mathematics, University of Arkansas,\\
Fayetteville, AR 72701, U. S. A.}
\maketitle
\begin{quote}
\begin{abstract}
\ Plemelj projection operators are introduced for spaces of square integrable
functions defined over the boundaries of a class of compact real n-dimensional  
manifolds lying in $C^{n}$. These manifolds posses many properties similar to
domains in $R^{n}$, and are consequently called domain manifolds. The
key ingredients used here are techniques from both real and complex Clifford
analysis. Analogues of the Kerzman-Stein kernel and Szeg\"{o} projection
operators are introduced, and their conformal covariance is described.
\end{abstract}
\end{quote}
\begin{center}{\Large Introduction}\end{center}

\ It is reasonably well known that given any reasonably smooth curve, $S$, 
lying in the complex plane and dividing $C$ into two open components,
then $L^{p}(S)=H^{p}(S^{+})\oplus H^{p}(S^{-})$ for $1<p<\infty$.
Here $S^{+}$ and $S^{-}$ are the respective domains 
that complement the curve $S$, and $H^{p}(S^{\pm})$ are the Hardy spaces of 
analytic functions defined on $S^{\pm}$ respectively and which extend 
continuously
in the $L^{p}$ sense to their common boundary, $S$. This decomposition of 
$L^{p}(S)$ is obtained using Plemelj projection operators, which in turn are 
described using the singular Cauchy transform over $S$ and Cauchy's integral 
formula, see for instance [B].

\ When turning to higher dimensions one can consider a sufficiently                                                 smooth, 
orientable surface, $S$, lying in $R^{n}$ and such that $S$ divides $R^{n}$
into two complementary domains, $S^{+}$ and $S^{-}$ respectively. Again one
may consider the space $L^{p}(S)$ of $p$-integrable functions over $S$ where 
$1<p<\infty$. By introducing the Clifford algebra generated from $R^{n}$, the 
Dirac operator over $R^{n}$ the associated Cauchy integral formula and 
singular Cauchy transform over $S$ one can obtain the decomposition
$L^{p}(S)=H^{p}(S^{+})\oplus H^{p}(S^{-})$. The spaces $H^{p}(S^{\pm})$ are 
the spaces of solutions to the Dirac equation, or generalized Cauchy-Riemann 
system, which are defined on $S^{\pm}$ respectively and continuously extend 
in the $L^{p}$ sense to their common boundary, $S$. This decomposition is 
outlined for the cases where $S$ is compact and either $C^{2}$ or Liapunov in 
[GuSp1], and for the cases where $S$ is a Lipschitz graph in [Mi, LMcQ]. 

\ In all these cases one needs to assume that $L^{p}(S)$ is the space of 
Clifford algebra valued $p$ integrable functions defined on $S$. As the 
Clifford algebra contains a unit then the usual complex space of $L^{p}$
is contained in this space.

\ This decomposition has proved to be a fundamental result in Clifford 
analysis and has provided a fundamental link between it and the study of
many types of boundary value problems. See for instance [Be, GuSp1, GuSp2, Mi,
LMcQ].

\ In greater generality one may use operator theory associated to general 
Dirac operators and associated Clifford algebras to obtain similar results 
over compact, orientable, smooth submanifolds of codimension $1$ lying
within a smooth, orientable, Riemannian manifold. See for instance [B-BWo, C]. 

\ In separate work, [R1] and references therein, the author has 
introduced a class of smooth, orientable,
real $n$-dimensional manifold which have many basic properties of domains 
lying in $R^{n}$, see also [GuKi, S]. Such manifolds are called domain manifolds.
In [GuKi] it is shown that one can set up and solve boundary value problems
over domain manifolds in much the same way as one does in complex analysis and 
Clifford analysis over Euclidean space. In particular one can set up Plemelj 
formulae, see [R1]. However, in [R1] these analogues of Plemelj formulae are
introduced only for spaces of H\"{o}lder continuous functions defined on the
boundary of a domain manifold, there referred to as a manifold of type one or 
two.

\ In this paper the intention is to introduce Plemelj projection operators 
and the analogous Hardy space decompositions over domain manifolds lying in 
both $C^{n}$. We also indicate how the results presented here extend to the 
analogues of domain manifolds lying in the complex sphere. Throughout we shall 
restrict attention to the cases where the boundary of the domain manifold is 
$C^{2}$.

\begin{center}{\Large Preliminaries}\end{center}

\ In this section we will develop the background on Clifford analysis that we 
shall need here.

\ We shall first consider the real, $2^{n}$ dimensional Clifford algebra, 
$Cl_{n}$, 
generated from $R^{n}$. So if $R^{n}$ has orthonormal basis 
$e_{1},\ldots, e_{n}$ then 
$Cl_{n}$ is chosen so that it has basis $1, e_{1},\ldots, e_{n},\ldots,
e_{j_{1}}\ldots e_{j_{r}},\ldots, e_{1}\ldots e_{n}$, where $1\leq j_{1}<
\ldots< j_{r}\leq n$ and $e_{i}e_{j}+e_{j}e_{i}=-2\delta_{i, j}$ where
$\delta_{i, j}$ is the Kroneker delta function. It may be seen that we 
are considering $R^{n}$ to be embedded in $Cl_{n}$. 

\ Although in general $Cl_{n}$ is not a division algebra it may be seen that
 each non-zero vector $x\in R^{n}$ has a multiplicative inverse $x^{-1}=
\frac{-x}{\|x\|^{2}}$. Up to a sign this inverse corresponds to the Kelvin 
inverse of a non-zero vector.

\ We shall need the following antiautomorphism over $Cl_{n}$ 
\[-:Cl_{n}\rightarrow Cl_{n}:-(e_{j_{1}}\ldots e_{j_{r}})=
(-1)^{r}e_{j_{r}}\ldots e_{j_{1}}.\]

\ For $a=a_{0}+\ldots+a_{1,\ldots, n}e_{1}\ldots e_{n}\in Cl_{n}$ we shall
write $\stackrel{-}{a}$ for $-a$. It
may be noted that the real part of $a\stackrel{-}{a}$ gives the square of 
the norm, $\|a\|$, of $a$, where $\|a\|=
(a_{0}^{2}+\ldots +a_{1,\ldots, n}^{2})^{\frac{1}{2}}$.

\ We shall also need the complexification $Cl_{n}(C)$ of $Cl_{n}$. The
antiautomorphism, $-$ extends linearly to this algebra. As is 
usual we shall denote the complexification of $R^{n}$ by $C^{n}$. It is no 
longer the case that each non-zero vector in $C^{n}$ has a multiplicative
inverse. For instance we may consider $e_{1}+ie_{2}$. In this case
$(e_{1}+ie_{2})^{2}=0$. We shall denote a vector 
$z_{1}e_{1}+\ldots+z_{n}e_{n}\in C^{n}$ by $z$, and we shall denote the null
cone $\{z\in C^{n}:z^{2}=0\}$ by $N(0)$. For $z_{1}\in C^{n}$ we denote the
null cone $\{z\in C^{n}:(z-z_{1})^{2}=0\}$ by $N(z_{1})$. It may be noted 
that each vector $z\in C^{n}\backslash N(0)$ has a multiplicative inverse.
\begin{definition}
A smooth real n-dimensional manifold $M$ lying in $C^{n}$ is called a domain 
manifold if for each $z\in M$ 
\newline
(i) $N(z)\cap M=\{z\}$
\newline
(ii) $N(z)\cap TM_{z}=\{z\}$.
\end{definition}

\ Such manifolds have been used extensively in [GuKi, R1, S]. Here we shall 
assume that the boundary of each domain manifold considered is at worst 
$C^{2}$, though similar results to the ones presented here can be obtained 
with only minor modification to the techniques used here if one assumed that 
the boundary is $C^{1}$ and has H\"{o}lder continuous derivative with respect 
to some atlas. We shall not consider the case where the boundary of $M$ is 
Lipschitz. When the domain manifold is a subset of $R^{n}$ then its 
interior is a domain in $R^{n}$. It is for this reason that such manifolds are 
called domain manifolds.

\ Associated to each domain manifold $M$ is a domain in $C^{n}$ called a cell
of harmonicity. This domain is the component of 
$C^{n}\backslash\cup_{z\in \partial M}N(z)$ which contains the interior of 
$M$. We shall denote this domain by $M^{\dagger}$. When $M$ is a domain in 
$R^{n}$ then $M^{\dagger}$ is the cell of harmonicity, or Vekua hull, 
described in [A] and elsewhere.
\begin{definition}
For $U$ a domain in $R^{n}$ a $Cl_{n}(C)$ valued differentiable function
$f(x)$ defined on $U$ is said to be left monogenic if $Df=0$ everywhere on 
$U$, where $D=\Sigma_{j=1}^{n}e_{j}\frac{\partial}{\partial x_{j}}$. Similarly
a $Cl_{n}(C)$ valued differentiable function $g$ defined on $U$ is said to
be right monogenic if $gD=0$ everywhere on $U$, where 
$gD=\Sigma_{j=1}^{n}\frac{\partial g}{\partial x_{j}}e_{j}$.
\end{definition}

\ The following two theorems are standard to Clifford analysis and can be 
found in [BDSo, GMu] and elsewhere. 
\begin{theorem}{\bf{(Cauchy's Theorem)}}
Suppose that $f$ and $g$ are respectively left and right monogenic functions 
on $U$ and that $S$ is a compact $C^{1}$ surface bounding a subdomain of $U$
then
\[\int_{S}g(x)n(x)f(x)d\sigma(x)=0\]
where $n(x)$ is the outward pointing normal vector to $S$ at $x$ and $\sigma$
is the usual Hausdorff measure on $S$.
\end{theorem}
\begin{theorem}{\bf{(Cauchy's Integral Formula)}}
Suppose that $f$ is a left monogenic function on the domain $U$ and that $S$ 
is a compact $C^{1}$ surface bounding a subdomain $V$ of $U$. Then for each 
$y\in V$
\[f(y)=\frac{1}{\omega_{n}}\int_{S}G(x-y)n(x)f(x)d\sigma(x)\]
where $\omega_{n}$ is the surface area of the unit sphere in $R^{n}$ and
$G(x)$ is the left and right monogenic function $\frac{x-y}{\|x-y\|^{n}}$.
\end{theorem}

\ If $U^{\dagger}$ is the cell of harmonicity associated to the domain $U$ then
the previous theorem can be applied to show that when $n$ is even the left 
monogenic function has a unique holomorphic continuation to a function 
$f^{\dagger}$ defined on $U^{\dagger}$. Moreover $f^{\dagger}$ satisfies the
the equation $D_{C}f^{\dagger}=0$ where 
$D_{C}=\Sigma_{j=1}^{n}e_{j}\frac{\partial}{\partial z_{j}}$. When $n$ is odd
then $f^{\dagger}$ is defined on some Riemann surface covering $U^{\dagger}$.
\begin{definition}
Suppose that $\Omega$ is a domain in $C^{n}$ or a Riemann surface covering a 
domain in $C^{n}$ and $f(z)$ is a holomorphic 
function defined on $\Omega$ and takes values in $Cl_{n}(C)$. Then if $f$ 
satisfies the equation $D_{C}f=0$ then $f$ is called a complex left monogenic 
function. Similarly if $g$ is a $Cl_{n}(C)$ valued holomorphic function
defined on $\Omega$ and satisfying $gD_{C}=0$ then $g$ is called a complex
right monogenic function. 
\end{definition}

\ The following two theorems are treated in depth in [R1, S].
\begin{theorem}
Suppose that $f$ and $g$ are respectively complex left and right monogenic
functions defined in a neighbourhood $\Omega$ of a compact domain manifold
$M$. Then
\begin{equation}
\int_{\partial M}g(z)Wzf(z)=0
\end{equation}
where $Wz$ is the holomorphic differential form 
$\Sigma_{j=1}^{n}(-1)^{n}e_{j}dz_{1}\wedge\ldots\wedge dz_{j-1}
\wedge dz_{j+1}\wedge\ldots\wedge dz_{n}$.
\end{theorem}
\begin{theorem}
Suppose that $n$ is even and that $f$ is a complex left monogenic function
defined in a neighbourhood of a compact domain manifold $M$. Then $f$ has a 
unique holomorphic continuation to a complex left holomorphic function 
$f^{\dagger}$ on $M^{\dagger}$ and for each $w\in M^{\dagger}$
\begin{equation}
f^{\dagger}(w)=\frac{1}{\omega_{n}}\int_{\partial M}G(z-w)Wzf(z)
\end{equation}
where $G(z)=(-1)^{\frac{n}{2}}z^{-n+1}$ is the complex left and right 
monogenic extension of $G(x)$ to $C^{n}\backslash N(0)$.
\end{theorem}

\ In the case where $n$ is odd the complex left monogenic function defined by 
Equation 2 is defined on some Riemann surface covering a domain 
containing $M$.

\ In [QR] it is pointed out that the differential form $Wz$ appearing in 
equations 1 and 2 can be replaced so that equation 1 becomes
\[\int_{\partial M}g(z)n(z)f(z)d\sigma(z)=0\]
while equation 2 becomes
\[f^{\dagger}(w)=\frac{1}{\omega_{n}}\int_{\partial M}g(w-z)n(z)f(z)
d\sigma(z)\] 
where $n(z)$ is a vector in $C^{n}$ orthogonal, with respect to the inner 
product $<z,w>=\Sigma_{j=1}^{n}z_{j}w_{j}$, to the complexification of the
tangent space $T\partial M_{z}$. Moreover $\sigma$ is a complex valued 
measure on $\partial M$. In [QR] it is observed that the function $n(z)$ is
non-zero and $C^{1}$. 

\begin{center}{\Large Hardy Spaces over Domain Manifolds}\end{center}

\begin{definition}
Suppose that $M$ is a domain manifold then for $1<p<\infty$ a function 
$f:\partial M\rightarrow Cl_{n}(C)$ is said to belong to the $L^{p}$ space, 
$L^{p}(\partial M)$, of $\partial M$ if 
\[(\int_{\partial M}\|f(z)\|^{p}|d\sigma(z)|)^{\frac{1}{p}}<\infty\].
\end{definition}

\ In the previous definition the term $|d\sigma(z)|$ stands for the 
infinitesimal of the absolute value measure $\sigma|$ of the complex measure 
$\sigma$ and it is introduced in [QR]. 

\ For each pair $f,g\in L^{2}(\partial M)$ there is a well defined
inner product
\[<f,g>_{\partial M}=\int_{\partial M}\overline{f(z)}g(z)d\sigma(z).\]
This integral is dominated by $\int_{\partial M}f(z)^{\star}g(z)d\sigma(z)$,
where $f^{\star}$ is the complex conjugate of $\overline{f}$. The identity 
component of 
\[\int_{\partial M}f^{\star}(z)f(z)|d\sigma(z)|\]
gives the
square of the $L^{2}$ norm of $f$ for each $f\in L^{2}(\partial M)$.
 
\ When $M$ is a subset of $R^{n}$ this definition of an $L^{p}$ space
corresponds to the usual definition of the $L^{p}$ space over the boundary
of some domain in $R^{n}$.

\ The following result is an immediate consequence of theorem 2.6 on page
223 of [StW].
\begin{proposition}
Suppose that $V$ is a real $n-1$ dimensional vector subspace of $C^{n}$. Then 
for $1<p<\infty$ the integral
\[P. V.\frac{1}{\omega_{n}}\int_{V}G(z-w)n(w)f(w)d\sigma(w)\]
defines a bounded linear operator
\[C_{V}:L^{p}(V)\rightarrow L^{p}(V).\]
\end{proposition}

\ Using this result we can now deduce:
\begin{theorem}
Suppose that $M$ is a compact domain manifold then the integral
\[P. V. \frac{1}{\omega_{n}}\int_{\partial M}G(z-w)n(w)f(w)d\sigma(w)\]
defines a bounded linear operator
\[C_{\partial M}:L^{2}(\partial M)\rightarrow L^{2}(\partial M).\]
\end{theorem}
{\bf{Proof:}} The proof follows a standard argument for compact, $C^{2}$, or 
even Liapunov, curves lying in the complex plane. This argument was adapted to 
the Clifford analysis setting for sufficiently smooth compact surfaces lying in 
$R^{n}$ in [Be]. The argument makes use of a cancellation property given by 
the $C^{1}$ function $n(w)$ over $\partial M$.

\ Essentially for each $z\in \partial M$ we may find an $\epsilon(z)\in R^{+}$
and a $C^{1}$ homotopy 
\[H_{z}:B(z,\epsilon(z))\cap\partial M\times[0,1]\rightarrow C^{n}\]
such that 
\newline
(i) $H(w,0)=w$
\newline
(ii) $H(B(z,\epsilon(z)\cap\partial M), 1)=
B(z,\epsilon(z))\cap T\partial M_{z}$.
\newline
(iii) $\|H(w,1)-z\|=\|w-z\|$.
\newline
Moreover,
\newline
(iv) there is a $\delta(z)\in R^{+}\cup\{0\}$ such that 
$\|w-H(w,1)\|<\delta(z)\|w-z\|$ for each $w\in B(z,\epsilon(z))\cap\partial M$.

\ Now 
\[\int_{B(z,\epsilon(z))\cap\partial M}G(z-w)n(w)f(w)d\sigma(w)=\]
\[\int_{B(z,\epsilon(z))\cap\partial M}G(z-w)(n(w)-n(z))f(w)d\sigma(w)\]
\[+\int_{B(z,\epsilon(z))\cap\partial M}G(z-w)n(z)f(w)d\sigma(w).\]

\ As $n(w)$ is a $C^{1}$ function over $\partial M$ then the first term on 
the right side of the previous equation defines a weakly singular integral
operator acting on the square integrable function $f$.

\ We can now use the homotopy $H$ and Proposition 1 to deal with the term
\[\int_{B(z,\epsilon(z))\cap \partial M}G(z-w)n(z)f(w)d\sigma(w)\]

\ This term can be rewritten as
\[\int_{B(z,\epsilon(z))\cap\partial M}(G(z-w)-G(z-H(w,1)))n(z)f(w)d\sigma(w)\]
\[+\int_{B(z,\epsilon(z))\cap\partial TM_{z}}
G(z-v)n(z)\lambda_{z}(v)f(\psi_{z}(v))d\sigma(v),\]
where $\psi_{z}$ is the $C^{1}$ diffeomorphism defined by $H(\psi_{z}(v),1)=v$, 
and $\lambda_{z}$ is the Jacobian associated to the $C^{1}$ function 
$\psi_{z}$.

\ By conditions (iii) and (iv) of the homotopy $H$ the first term on the right 
side of the previous equation is dominated by the term
\[C\delta(z)\int_{B_(z,\epsilon(z))\cap\partial M}\|z-w\|^{-n+2}\|n(z)f(w)\|
|d\sigma(w)|\]
for some dimensional constant $C\in R^{+}$, and so again defines a weakly 
singular integral operator acting on $f$. 

\ As $\lambda_{z}$ is a bounded measurable function it follows from
proposition 1 that the second term on the right side of the same equation also
defines an $L^{2}$ bounded operator. 

\ Moreover, the term
\[\int_{\partial M\backslash(B(z,\epsilon(z))\cap\partial M)}G(z-w)n(w)f(w)
d\sigma(w)\]
is $L^{2}$ bounded.

\ So far we have shown that for each $z\in\partial M$ the operator
\[T_{z}:L^{2}(\partial M)\rightarrow L^{2}(\partial M)\]
defined by
\[\int_{\partial M\cap B(z,\epsilon(z))}G(z-w)n(w)f(w)d\sigma(w)\]
for each $f\in L^{2}(\partial M)$, is an $L^{2}$ bounded operator. The result 
now follows from the compactness of $\partial M$. $\Box$

\ If in Proposition 1 and Theorem 5 we replace the kernel $G(z)$ by a kernel 
$K(z)=\frac{\Omega(z)}{(z^{2})^{\frac{n-2}{2}}}$ where $\Omega(z)$ is an odd 
function homogeneous of degree zero, we also obtain an $L^{2}$ bounded 
operator
\[T_{K}:L^{2}(\partial M)\rightarrow L^{2}(\partial M):
T_{K}(g)(w)=P. V. \int_{\partial M}K(w-z)n(z)g(z)d\sigma(z).\]
The proof is the same as the proof of Theorem 5.

\ Suppose that $g$ belongs to the function space $C^{1}(\partial M)$
of $Cl_{n}(C)$ valued $C^{1}$ functions on $\partial M$. Suppose also that
 $\theta:[0, 1)\rightarrow M$ is piecewise smooth with 
$lim_{t\rightarrow 1}\theta(t)=w\in\partial M$ and $\theta'(1)$ does not
belong to $TM_{w}(C)$, 
the complexification of the tangent space $TM_{w}$. Then in [R1] we show that
\[\lim_{t\rightarrow 1}\frac{1}{\omega_{n}}\int_{\partial M}
G(\theta(t)-z)n(z)g(z)d\sigma(z)=\]
\[\frac{1}{2}g(w)
+\frac{1}{\omega_{n}} P. V. \int_{\partial M}G(w-z)n(z)g(z)d\sigma(z).\]

\ Using the fact that $C^{1}(\partial M)$ is a dense subspace of 
$L^{2}(\partial M)$ we immediately have:
\begin{theorem}
Suppose that $f\in L^{2}(\partial M)$ and that $\theta$ and $w$ are as in 
the preceding paragraph. Then
\[\lim_{t\rightarrow 1}\frac{1}{\omega_{n}}\int_{\partial M}
G(\theta(t)-z)n(z)f(z)d\sigma(z)=\]
\[\frac{1}{2}f(w)+\frac{1}{\omega_{n}}P. V. 
\int_{\partial M}G(w-z)n(z)f(z)d\sigma(z) \hspace{.5in}a.e..\]
\end{theorem}

\ This last formula is a generalization of the Plemelj formula arising in
one variable complex analysis.

\ We would like to show that the pointwise convergence described in the
previous equation can be replaced by uniform convergence in the $L^{2}$ norm.
To this end we first introduce maximal functions.
\begin{definition}
Suppose that $f\in L^{p}(\partial M)$ and $1<p<\infty$ then for each 
$w\in\partial M$ we define the maximal function $M(f)(w)$ of $f(w)$to be 
\[\sup_{r>0}\frac{1}{\int_{B(w, r)\cap\partial M}|d\sigma(z)|}
\int_{B(w, r)\cap\partial M}\|f(z)\||d\sigma(z)|.\]
\end{definition}

\ By similar arguments to those presented in [St] it may be deduced that for
each $f\in L^{p}(\partial M)$ and with $1<p<\infty$ then 
$M(f)\in L^{p}(\partial M)$ and there is a positive constant $C(n, p)$ such 
that
\begin{equation} 
\|M(f)\|_{L^{p}}\leq C(n, p)\|f\|_{L^{p}}.
\end{equation} 

\ We also will need non-tangential maximal functions. By identifying $C^{n}$
with $R^{2n}$ endowed with the usual inner product $<,>_{R^{2n}}$ for 
each $w\in C^{n}$ each $u\in C^{n}\backslash\{0\}$ and each 
$\alpha\in (0,\frac{\pi}{2})$ we can introduce the cone 
$\Gamma(w,u, \alpha)=\{z\in C^{n}:0<<z-w,u>_{R^{2n}}<\|z-w\|\cos\alpha\}$.
For each $r\in R^{+}$ we define the truncated cone $\Gamma(w, u,\alpha, r)$
to be the set $\{z\in \Gamma(w, u, \alpha):\|z-w\|<r\}$. Both the sets
$\Gamma(w, u,\alpha)$ and $\Gamma(w, u, \alpha, r)$ are open subsets of 
$C^{n}$. 

\ As the manifold $M$ is compact and has $C^{2}$ boundary then we may find
an $\alpha(M)\in (0, \frac{\pi}{2})$ and an $r(M)\in R^{+}$ such that
$\Gamma(w, n(w), \alpha(M), r(M))\subset M^{\dagger}$ and
$\partial\Gamma(w, n(w),\alpha(M), r(M))\cap\partial M^{\dagger}=\{w\}$.
for each $w\in\partial M$.
\begin{definition}
For each $p\in(0,\infty)$ and each $f\in L^{p}(\partial M)$ we define the 
non-tangential maximal function $N(f)(w)$ of $f$ to be 
\[\sup_{z\in\Gamma(w, n(w),\alpha(M), r(M))}\|
\int_{\partial M}G(z-u)n(u)f(u)d\sigma(u)\|\]
\end{definition}

\ The non-tangential maximal function introduced in the previous definition 
differs from the usual one set up over domains in $R^{n}$. There one 
considers supremums over real n-dimensional cones lying in $R^{n}$. Here we 
have doubled the dimension and are taking supremums over complex, 
n-dimensional cones in $C^{n}$.

\begin{theorem}
Suppose that $f\in L^{2}(\partial M)$ then 
$\|N(f)\|_{L^{2}}<C(M)\|f\|_{L^{2}}$ for some $C(M)\in R^{+}$.
\end{theorem}
{\bf{Proof:}} The proof is an adaptation of arguments appearing on page 63 of
[Mi] and page 27 of [MCo].

\ Let us begin by considering 
\[\int_{\partial M\backslash (B(w, \epsilon)\cap\partial M}
G(w-z)n(z)f(z)d\sigma(z)-\int_{\partial M}G(u-z)n(z)f(z)d\sigma(z),\]
where $w\in\partial M$ and $u\in M^{\dagger}$.
This expression is equal to
\[\int_{\partial M\backslash (B(w,\epsilon)\cap\partial M}
(G(w-z)-G(u-z))n(z)f(z)d\sigma(z)\]
\[-\int_{B(w,\epsilon)\cap\partial M}G(u-z)n(z)f(z)d\sigma(z).\]
Let us restrict $u$ so that $\|w-u\|=\epsilon$ and 
$u\in\Gamma(w,n(w),\alpha(M), r(M))$.
Now
\[\|\int_{B(w,\epsilon)\cap\partial M}G(u-z)n(z)f(z)d\sigma(z)\|\leq
C\frac{1}{\epsilon^{n-1}}\int_{B(w,\epsilon)\cap\partial M}
\|f(z)\||d\sigma(z)|\]
\[\leq C_{1}M(f)(w),\]
for some $C$ and $C_{1}\in R^{+}$.

\ Moreover
\[\|G(w-z)-G(u-z)\|\leq C_{2}\frac{\epsilon}{\|z-w\|^{n}}\]
for each $z\in\partial M\backslash(B(w,\epsilon)\cap\partial M$ and some 
constant $C_{2}\in R^{+}$. 

\ It follows that
\[\|\int_{\partial M\backslash (B(w,\epsilon)\cap\partial M)}(G(w-z)-G(u-z))
n(z)f(z)d\sigma(z)\|\leq\]
\[C_{2}\epsilon\sum_{j=0}^{\infty}
\|\int_{\partial M\cap(B(w,2^{j}\epsilon)\backslash B(w,2^{j+1}\epsilon))}
\frac{\|f(z)\|}{\|z-w\|^{n}}|d\sigma(z)|,\]
for some $C_{2}\in R^{+}$. The right side of this expression is dominated by
$C_{3}M(f)(w)$ for some $C_{3}\in R^{+}$.

\ Consequently
\[N(f)(w)<C_{4}(M(f)(w)
+\sup_{\epsilon>0}\|\int_{B(w,\epsilon)\cap\partial M}
G(w-z)n(z)f(z)d\sigma(z)\|),\]
for some $C_{4}\in R^{+}$. This inequality is derived in the euclidean setting
using much the same arguments in [Mi].

\ Now
\[\|\int_{\partial M\backslash(B(w,\epsilon)\cap\partial M)}
G(w-z)n(z)f(z)d\sigma(z)\|
\leq(\|\int_{\partial M}G(w-z)n(z)f(z)d\sigma(z)\|\]
\[+\|\int_{B(w,\epsilon)\cap\partial M}G(w-z)n(z)f(z)d\sigma(z)\|.\]

\ One can now readily adapt the proof of Cotlar's inequality given in [MCo]
and show that
\[\sup_{\epsilon>0}\|\int_{\partial M\backslash B(w,\epsilon)}G(w-z)n(z)
f(z)d\sigma(z)\|\leq \]
\[C(M(P. V. \int_{\partial M}G(w-z)n(z)f(z)d\sigma(z))+M(f)(w))\]
for some constant $C\in R^{+}$.

\ The result now follows from Theorem 5 and Inequality 3. $\Box$

\ Combining Theorems 6 and 7 with Lebesgue's dominated convergence theorem
we get:
\begin{theorem}
Suppose that $f\in L^{2}(\partial M)$ then for any smooth homotopy 
deformation $H:\partial M\times[0,1)\rightarrow M^{\dagger}$ such that
$\lim_{t\rightarrow 1}H(w,t)=w$ for each $w\in \partial M$
\[\lim_{t\rightarrow 1}\|\frac{1}{\omega_{n}}
\int_{\partial M}G(H(w,t)-z)n(z)f(z)d\sigma(z)\]
\[-\frac{1}{2}f(w)-\frac{1}{\omega_{n}}P. V. 
\int_{\partial M}G(w-z)n(z)f(z)\|_{L^{2}}=0.\]
\end{theorem}

\ Even in the special cases where $M$ is an open subset $U$ of $R^{n}$ 
Theorem 8 is a stronger statement than the usual statement of uniform 
convergence in the $L^{2}$ norm. This is because one usually looks for 
uniform convergence within the domain $U$ and not on the cell of harmonicity 
$U^{\dagger}$ associated to $U$.

\ Let us denote the $L^{2}$ bounded operator defined by
\[\frac{1}{2}f(w)+\frac{1}{\omega_{n}}
\int_{\partial M}G(w-z)n(z)f(z)d\sigma(z)\]
by $S^{+}$ and the singular integral operator $S^{+}-\frac{1}{2}I$ by 
$C_{\partial M}$. For each $f\in L^{2}(\partial M)$ let us denote $S^{+}f$ by
$f^{+}$. A simple application of Cauchy's integral formula and Theorem 6
reveal that $S^{+2}=S^{+}$ and $C_{\partial M}^{2}=\frac{1}{4}I$. If we denote 
the $L^{2}$ bounded operator $\frac{1}{2}I-C_{M}$ by $S^{-}$ then we
can readily deduce that $S^{-2}=S^{-}$ and $S^{+}S^{-}=S^{-}S^{+}=0$. This 
generalizes to $L^{2}(\partial M)$ results known to hold over the $L^{2}$ 
space of a closed reasonably smooth curve in the complex plane and for $L^{2}$
spaces of sufficiently smooth surfaces lying in $R^{n}$. See for instance [B,
LMcQ] and elsewhere.

\ In [R1] we note that for any compact domain manifold $M$ we can construct
another domain manifold $M'$ such that $M\subset M'$ and $\partial M$ is 
contained in the interior of $M'$. Let us suppose that $\partial M$ is
connected then let us denote the component of 
$C^{n}\backslash\cup_{z\in\partial M}N(z)$ containing $M'\backslash cl(M)$ by
$M^{\dagger'}$, where $cl(M)$ is the closure of $M$. By arguments described in 
[R1] it may be noted that $M^{\dagger'}$ is a domain in $C^{n}$. In the cases
where $\partial M$ is multiply connected then $M^{\dagger'}$ is a union of 
finitely many disjoint domains.

\ By similar arguments to those used to deduce Theorem 6 one may deduce:
\begin{theorem}
Suppose that $f\in L^{2}(\partial M)$ then
\[\lim_{u\rightarrow w}\frac{1}{\omega_{n}}\int_{\partial M}G(u-z)n(z)f(z)
d\sigma(z)=-S^{-}f(w)\hspace{0.5in}a.e.\]
where $u\in M^{\dagger'}$, $w\in\partial M$ and $u$ approaches $w$ 
non-tangentially.
\end{theorem}

\ If we denote $-S^{-}f$ by $f^{-}$ for each $f\in L^{2}(\partial M)$, then
$f=f^{+}+f^{-}$.

\ By similar arguments to those used to deduce Theorem 8 we can also derive:
\begin{theorem}
Suppose that $H:\partial M\times [0,1)\rightarrow M^{\dagger'}$ is a smooth 
homotopy such that $\lim_{t\rightarrow 1}H(w,t)=w$ then
\[\lim_{t\rightarrow 1}\|\frac{1}{\omega_{n}}
\int_{\partial M}G(H(w,t)-z)n(z)f(z)d\sigma(z)-f^{-}(w)\|_{L^{2}}=0\]
for each $f\in L^{2}(\partial M)$.
\end{theorem}
\begin{definition}
We shall call the space $S^{+}L^{2}(\partial M)$ the Hardy 2-space of 
$M^{\dagger}$ and we shall denote it by $H^{2}(M^{\dagger})$. 

\ Similarly we
shall call the space $S^{-}L^{2}(\partial M)$ the Hardy 2-space of 
$M^{\dagger'}$. We shall denote this space by $H^{2}(M^{\dagger'})$.
\end{definition}

\begin{theorem}
\[L^{2}(\partial M)=H^{2}(M^{\dagger})\oplus H^{2}(M^{\dagger'}).\]
\end{theorem}

\ Let $P^{+}:L^{2}(\partial M)\rightarrow H^{2}(M^{\dagger})$ be the 
Szeg\"{o}, or orthogonal, projection with respect to the inner product
$<,>_{\partial M}$. Also let 
$P^{-}:L^{2}(\partial M)\rightarrow H^{2}(M^{\dagger'})$ be the Szeg\"{o}, or 
orthogonal, projection with respect to the same inner product 
$<,>_{\partial M}$. These operators are self adjoint with respect to this
inner product. In other words $P^{\pm\star}=P^{\pm}$ where $P^{\pm\star}$ is 
the adjoint of $P^{\pm}$ with respect to $<,>_{\partial M}$.

\ Also let $C_{\partial M}^{\star}(f)(w)=
\frac{1}{\omega_{n}}n(w)\int_{\partial M}G(w-z)f(z)d\sigma(z)$. Then the 
$L^{2}$ bounded operator $C_{\partial M}-C_{\partial M}^{\star}$ is called 
the Kerzman-Stein kernel. This is in complete analogy to the complex variable 
setting, see [B].

\ We may easily obtain the following identities.
\newline
$P^{\pm}=S^{\pm}P^{\pm}$, $P^{\pm\star}=P^{\pm\star}S^{\pm\star}=
P^{\pm}S^{\pm\star}$, and $S^{\pm}=P^{\pm}S^{\pm}$. 
Consequently $P^{\pm}-S^{\pm}=
P^{\pm}(C_{\partial M}^{\star}-C_{\partial M})$.

\ Hence
\[P^{\pm}(I-(C_{\partial M}^{\star}-C_{\partial M}))=S^{\pm}.\]
These identities are generalizations of identities appearing in [B, Se].

\ As in the classical case, see [KSt], the singularities of the integral
operators $C_{\partial M}$ and $C_{\partial M}^{\star}$ are cancelled out
in $C_{\partial M}-C_{\partial M}^{\star}$. Consequently this Kerzman Stein 
kernel is a compact operator.

\ In order to move from the setting where $\partial M$ is compact to the more
general setting we shall use M\"{o}bius transformations. In [Bo] and elswhere 
it is shown that the space of left monogenic functions over some domain remain 
invariant under M\"{o}bius transformations. Also in [QR] it is noted that 
cells of harmonicity transform to other cells of harmonicity under M\"{o}bius 
transformations. So in particular if $u=\psi(z)=(z+a)^{-1}$ for some constant 
$a\in C^{n}$ and $f(u)$ is a complex left monogenic function defined on a
cell of harmonicity $M^{\dagger}$ then $G(z+a)f((z+a)^{-1})$ is a complex left
monogenic function on the cell of harmonicity $\psi^{-1}(M^{\dagger})$. Also, 
[QR], $g(u)\in L^{2}(\partial M)$ if and only if 
$G(z+a)g((z+a)^{-1})\in L^{2}(\partial\psi^{-1}(M))$. Moreover, [QR], 
$\|G(z+a)g((z+a)^{-1})\|_{L^{2}}=\|g(u)\|_{L^{2}}$. It follows that
$\|M(G(z+a)g((z+a)^{-1}))\|_{L^{2}}\leq\|g(u)\|_{L^{2}}$ for each 
$g\in L^{2}(\partial M)$.

\ In the case where $-a\in\partial M$ then $\partial M$ is no
longer compact. Even though $\partial M$ is no longer compact its $C^{2}$
structure is preserved.

\ In [QR] it is noted that if $v=(w+a)^{-1}$ then $G(v-u)=
G(w)^{-1}G(w-z)G(z)^{-1}$ and if $f$, $g\in L^{2}(\partial M)$ then
\[\int_{\partial M}f(u)n(u)g(u)d\sigma(u)=\]
\[\int_{\partial\psi^{-1}(M)}
f((z+a)^{-1})G(z+a)n(z)G(z+a)g((z+a)^{-1})d\sigma(z).\]

\ From these remarks it is an easy exercise to transpose all results so far 
obtained here for the $L^{2}$ space of a compact $C^{2}$ manifold $\partial M$
to the setting where $\partial M$ is no longer compact but is the image under
a M\"{o}bius transformation $\psi$ of a compact $C^{2}$ boundary of a domain 
manifold.

\ In [R2] we use a Cayley transformation to show that much of Clifford analysis
over $R^{n}$ can also be set up over the sphere $S^{n}$ lying in $R^{n+1}$. 
Using this transformation it is a reasonably easy exercise to transpose the
results presented here over domain manifolds and cells of harmonicity to their
images under this Cayley transformation within the complex sphere
$S_{C}^{n}=\{z\in C^{n+1}:z^{2}=-1\}$.

\begin{center}{\Large References}\end{center}
[A] V. Avanissian, {\em{Cellule D'Harmonicit\'{e} et Prolongement 
Analytique Complexe}}, Hermann, Paris, 1985. 
\newline
[B] S. Bell, {\em{The Cauchy Transform, Potential Theory and Conformal
Mapping}}, CRC Press, Boca Raton, 1992. 
\newline
[Be] S. Bernstein, {\em{The left-linear Riemann problem in Clifford analysis}},
Bull. Belg. Math. Soc., 2, 1996, 557-576.
\newline
[Bo] B. Bojarski, {\em{Conformally covariant differential operators}}, 
Proceedings, XXth Iranian Math. Congress, Tehran, 1989.
\newline
[B-BWo] B. Booss-Bavnbeck and K. Wojciechowski {\em{Elliptic Boundary Problems
for Dirac Operators}}, Birkha\"{u}ser, Basel, 1993.
\newline
[BDSo] F. Brackx, R. Delanghe and F. Sommen, {\em{Clifford Analysis}}, Pitman 
Research Notes in Mathematics, No, 76, London, 1982.
\newline
[C] D. Calderbank, {\em{Clifford analysis for Dirac operators of manifolds 
with boundary}}, to appear.
\newline
[GMo] J. Gilbert and M. A. M. Murray, {\em{Clifford Algebras and Dirac 
Operators in Harmonic Analysis}}, CUP, Cambridge, 1991.
\newline
[GuKi] K. G\"{u}rlebeck and F. Kippig, {\em{Complex Clifford analysis and 
elliptic boundary value problems}}, Advances in Applied Clifford Algebra, 
5, 1995, 51-62.
\newline
[GuSp1] K. G\"{u}rlebeck and W. Spr\"{o}ssig, {\em{Quaternionic Analysis and 
Elliptic Boundary Value Problems}}, Birkha\"{u}ser Verlag, Basel, 1990.
\newline
[GuSp2] K. G\"{u}rlebeck and W. Spr\"{o}ssig, {\em{Quaternionic and Clifford
Calculus for Physicists and Engineers}}, Wiley and Sons, New York, 1997.
\newline
[LMcQ] C. Li, A. McIntosh and T. Qian, {\em{Clifford algebras, Fourier 
transforms and singular convolution operators on Lipschitz surfaces}}, 
Rev. Mat. Iberoamericana, 10, 1994, 665-721.
\newline
[KSt] N. Kerzman and E. M. Stein, {\em{The Cauchy kernel, the Szeg\"{o}
kernel, and the Riemann mapping function}}, Math. Ann., 236, 1978, 85-93.
\newline
[MCo] Y. Meyer and R. Coifman, {\em{Wavelets, Calder\'{o}n-Zygmund and
Multilinear Operators}}, CUP, Cambridge, 1997.
\newline
[Mi] M. Mitrea, {\em{Singular Integrals, Hardy Spaces, and Clifford 
Wavelets}}, Lecture Notes in Mathematics, No 1575, Springer Verlag, 
Heidelberg, 1994.
\newline
[QR] T. Qian and J. Ryan, {\em{Conformal transformations and Hardy spaces 
arising in Clifford analysis}}, Journal of Operator Theory, 35, 1996, 
349-372.
\newline
[R1] J. Ryan, {\em{Plemelj formulae and transformations associated to plane 
wave decompositions in complex Clifford analysis}}, Proceedings of the London 
\newline
Mathematical Society, 64, 1992, 70-94.
\newline
[R2] J. Ryan, {\em{Dirac operators on spheres and hyperbolae}}, Bolletin de la
Sociedad Matematica a Mexicana, 3, 1996, 255-270.
\newline
[S] K. Sano, {\em{Another type of Cauchy's integral formula in complex
Clifford analysis}}, Tokyo Journal of Mathematics, 20, 1997, 187-204.
\newline
[Se] S. Semmes, {\em{Chord-arc surfaces with small constant, 1}}, Advances 
in Mathematics, 85, 1991, 198-223.
\newline
[StW] E. M. Stein and G. Weiss, {\em{Introduction to Fourier Analysis on 
Euclidean Space}}, Princeton University Press, Princeton, 1971.
\newline
[St] E. M. Stein, {\em{Singular Integrals and Differentiability Properties of
Functions}}, Princeton University Press, Princeton, 1970.
\end{document}